# ON ADAPTIVE CONFIDENCE ELLIPSOIDS FOR SPARSE HIGH DIMENSIONAL LINEAR MODELS

**Xiaoyang Xie**
Department of Pure Mathematics and Mathematical Statistics
University of Cambridge
UK
xx765@cam.ac.uk

October 24, 2023

## ABSTRACT

In high-dimensional linear models the problem of constructing adaptive confidence sets for the full parameter is known to be generally impossible. We propose re-weighted loss functions under which constructing fully adaptive confidence sets for the parameter is shown to be possible. We give necessary and sufficient conditions on the weights for adaptive confidence sets to exist, and exhibit a concrete rate-optimal procedure in the feasible regime.

## 1 Introduction

We consider the high dimensional linear model

$$Y = X\theta + \varepsilon \tag{1}$$

where $Y \in \mathbb{R}^n$, $\theta \in \mathbb{R}^p$ and $X$ is a $n \times p$ matrix with entries $X_{ij}$. The interest is in statistically recovering some aspect of the unknown true parameter $\theta \in \mathbb{R}^p$ that generated the equation (1). In high dimensional problems where $n \to \infty$ and potentially $p \gg n$, it is common to assume the existence of some lower dimensional structure to achieve non-trivial results. We consider sparsity and assume that $\theta$ has $k$ non-zero coordinates:

$$\|\theta\|_0 \equiv \sum_{j=1}^{p} \mathbf{1}(\theta_j \neq 0) = k, k \in \mathbb{N},$$

for some $1 \leq k \leq p$. This is commonly referred to as $\theta$ being 'k-sparse' or $\theta$ being in the $l_0$-'ball' $B_0(k) := \{\theta : \|\theta\|_0 \leq k\}$.

In this paper we study the task of doing inference under the minimax paradigm, and we assume the following: $X_{ij} \overset{\text{i.i.d.}}{\sim} N(0,1)$, for $1 \leq i \leq n$, $1 \leq j \leq p$; and the noise $\varepsilon$ is isotropic standard Gaussian independent of $X$. Also assume that $n, p, k \to \infty$ with $k = o(p)$ and $p \simeq n$. The assumptions are made to focus on the exposition of the main statistical ideas. More involved setups, such as correlated design, unknown noise level or $p \lesssim \exp n$ can all be addressed in principle within our framework.

Assuming sparsity of $\theta$, under $\ell_2$-loss the estimation of the parameter is possible and the minimax optimal estimation rate is known to be $\sqrt{n^{-1}k\log p}$. Perhaps more importantly, plenty of estimators are adaptive to unknown sparsity: if $\theta \in B_0(k_{max})$ for some maximum sparsity level $k_{max} = o(n/\log p)$, the estimators achieve estimation rate $\sqrt{n^{-1}k\log p}$ for all $k < k_{max}$ without requiring the knowledge of the true sparsity $k$, see for example [1] and the in-depth analysis in [2].

We are interested in statistical confidence sets that adapt to the sparse estimation rate: ideally the diameter of the confidence sets matches said estimation rate without prior knowledge of $k$. To be precise: $C_n \equiv C_{n,\beta,\beta'}$, a measurable random subset of the parameter space $\Theta$ is called an *honest* confidence set if it satisfies

$$\inf_{\theta \in \Theta} \mathbb{P}_\theta(\theta \in C_n) \geq 1 - \beta$$

for any $\beta > 0$, while for every $\beta' > 0$ there exists a fixed constant $D$ s.t.

$$\limsup_{n\to\infty} \sup_{\theta\in\Theta} \mathbb{P}_\theta\Big(|C_n| > Dr_n(k)\Big) \leq \beta', 1 \leq k \leq k_{max}.$$

Here $r_n(k)$ is the estimation rate over $B_0(k)$. Constructing such a $C$ is known to be much more difficult than adaptive estimation and will be the focus of our study.

As shown and discussed in [3], the fundamental difficulty of uncertainty quantification under the sparse linear model arises from not being able to test between sub-models of different sparsity levels:

$$H_0 : \theta \in B_0(k_1) \quad vs \quad H_1 : \theta \in B_0(k_2) \setminus B_0(k_1)$$

for some $k_1 = o(k_2)$. Some separation between the hypotheses is necessary for effective and consistent testing, [4] developed the techniques for quantifying the separation and proving the indistinguishability result under insufficient separation. [5] further showed that the separation is closely related to the diameter of adaptive confidence sets: a consistent test coupled with an adaptive estimator leads to the existence of honest adaptive confidence sets. See Chapter 8.3 in [6] for comprehensive discussion of the relationship between minimax testing and uncertainty quantification.

The difficulty of constructing adaptive confidence sets reveals the 'geometric complexity' introduced by the sparsity assumption and its interaction with the loss function one chooses to do inference in. Previous work adaptive confidence sets have focused on fundamentally important choices of loss functions, namely under the $\ell_2$-loss

$$\|\vartheta\|_2^2 \coloneqq \|\vartheta\|_{\ell_2}^2 = \sum_{j=1}^{p} \vartheta_j^2, \vartheta \in \mathbb{R}^p \tag{2}$$

and the $l_\infty$-loss

$$\|\vartheta\|_\infty \coloneqq \max_{1\leq j\leq p} |\vartheta_j|, \vartheta \in \mathbb{R}^p. \tag{3}$$

For $\ell_\infty$-loss, or similarly any pre-specified $\theta_j$, the minimax estimation rate is almost the parametric rate (with a logarithmic penalty) and honest confidence sets 'adaptive' to the estimation rate (which does not depend on $k$) can be constructed. This can be achieved by, for example, constructing a class of de-sparsified estimators, as discussed in [7],[8] and [9]. However, under $\ell_2$-loss where the minimax estimation rate is $\sqrt{n^{-1}k\log p}$, [3] showed that confidence sets which adapt to rates faster than $p^{1/4}n^{-1/2}$ cannot exist (when $p \simeq n$).

Intuitively speaking, in the presence of the sparsity assumption and under common quadratic loss functions like the $\ell_2$-loss, the 'oracle' estimation rate becomes very fast relative to the dimensionality of the problem while the testing rate cannot be accelerated by the same amount unless certain part of the parameter space is removed. This is often referred to as 'separation', as $H_0$ and $H_1$ need to be sufficiently apart for testing to be possible. On the other hand, the de-sparsified estimators generate unbiased estimation of each coordinate that possesses good uncertainty quantification as a result of asymptotic normality, at the cost of *not* utilizing sparsity for faster estimation. Previously, adaptivity in uncertainty quantification was considered impossible, except in cases where the estimation rate collapses into the (almost) parametric rate: for example when doing inference on a single coordinate $\theta_j$ or on an 'ultra-sparse' linear functional as discussed in [10].

In some sense the $l_\infty$ case represents a task of doing inference on a small part of the parameter, while the $\ell_2$ case stands for inference on the full vector $\theta$. In this paper we study how the *transition* between the two scenarios occurs. For $\vartheta \in \mathbb{R}^p$, we take the loss function

$$\|\vartheta\|_{\ell_2^{-\alpha}}^2 \coloneqq \sum_{j=1}^{p} (j^{-\alpha}\vartheta_j)^2, \alpha \geq 0. \tag{4}$$

In Theorem 2.1 we prove that the adaptive estimation rate in $\ell_2^{-\alpha}$ is $\sqrt{n^{-1}k^{(1-2\alpha)}\log p}$ for $0 \leq \alpha \leq 1/2$ and $\sqrt{n^{-1}\log p}$ for $\alpha > 1/2$. Moving onto uncertainty quantification, it is shown in Theorem 2.2 that honest adaptive confidence sets can be constructed for all parameters if $\alpha \geq 1/4$. For $\alpha > \frac{1}{2}$, adaptation is trivial since estimation rate is not affected by true sparsity. When $1/4 \leq \alpha \leq 1/2$, our construction of honest adaptive confidence sets is based on risk estimation via the $U$-statistic, following the ideas in [11], [12], [3] and [13]. In the present context the $U$-statistic $U_n^{-\alpha}(\cdot)$ is defined as follows: For $\vartheta \in \mathbb{R}^p$, let

$$U_n^{-\alpha}(\vartheta) \coloneqq \frac{2}{n(n-1)} \sum_{j=1}^{p} \sum_{l=2}^{n} \sum_{m=1}^{l-1} j^{-2\alpha}(Y_l X_{lj} - \vartheta_j)(Y_m X_{mj} - \vartheta_j). \tag{5}$$

Meanwhile for any $0 < \alpha < 1/4$, we prove that adaptive confidence sets adapting to all possible sparsity levels cannot exist, due to obstacles similar in nature to under the $\ell_2$ case shown in [3].

As mentioned above, the existence of, or lack thereof, honest adaptive confidence sets depends on testing between sparsity levels. Under the loss functions (4) proposed here and for $0 \leq \alpha \leq \frac{1}{2}$, the estimation rate is $\sqrt{n^{-1}k^{1-2\alpha}\log p}$ while the 'critical' separation required for consistent testing is $p^{\frac{1}{4}-\alpha}n^{-1/2}$. For $\alpha < \frac{1}{4}$, the estimation rate is always faster than this rate in highly sparse sub-models so adaptation over all possible sparsity levels cannot be achieved. However, the testing rate accelerates relative to the estimation rate as $\alpha$ grows, creating a region where real adaptive uncertainty quantification is possible: when $\alpha \geq \frac{1}{4}$, rate-optimal confidence set exists at a genuine sparse-estimation rate without having to remove *any* part of the parameter space, and the phase transition occurs at precisely $\alpha = \frac{1}{4}$, as we show. It is conceivable that for $0 < \alpha < \frac{1}{4}$ partial adaptation over some range of $k$ can be achieved but we do not pursue this in the current work, as further development of specifying ranges of adaptation is mainly technical.

Seemingly, the $\ell_2^{-\alpha}$-loss down-weights each coordinate by a magnitude depending on how one chooses to order the $\theta_j$'s. It would be more pragmatic if information in the most significant coordinates is preserved. In practical problems it is very common that the interest lies in statistically recovering the top $m\%$ factors that carry the strongest signal for some pre-specified $m$. The positions of the important factors will be unknown, so ideally one wants to do inference under the loss functions

$$\sum_{j=1}^{p}(j^{-\alpha}\theta_{(j)})^2,$$

where $\theta_{(j)}$ are the *re-ordered* coordinates such that $|\theta_{(1)}| \geq |\theta_{(2)}| \geq \cdots \geq |\theta_{(p)}|$. The value and ordering of $\theta_j$'s are not directly observable but we show that it is indeed possible to perform inference in the *ordered* $\ell_2^{-\alpha}$-loss, as formally laid out in Corollary 2.3. Intuitively, this is expected since in proving the main theorems we have already shown that the performance of the inference procedures achieves the desired rate under the least favorable un-ordered case, i.e. where the non-zero entries of $\theta$ all happen to be in the first few coordinates.

Some of our ideas are related to [13] who study uncertainty quantification in density estimation, where the authors presented results of a related nature. The expression (4) is a Sobolev-style norm with a negative smoothness index, which was related to Wasserstein distances adopted in the aforementioned research. Of course, the geometric complexity arising from sparsity is distinct but the formal analogy remains.

Other than constructing adaptive confidence 'balls' for the vector $\theta$, a very closely related topic is constructing confidence intervals for linear functionals of $\theta$, namely for $L^T\theta$, $L \in \mathbb{R}^p$. It is possible to have fully adaptive confidence intervals for $L^T\theta$ when $L$ is very sparse, see [10] for example, but for general $L$ it is also known to be difficult. Consider a prediction problem where the interest is in random linear functionals: for example, as proposed in [14], we want to do inference on $X_{n+1}\theta$, $X_{n+1}^T \in \mathbb{R}^p$ with independent $N(0,1)$ entries. Realizing that $\mathbb{E}[(X_{n+1}\theta)^2] = \|\theta\|_2^2$, one can show that having a fully adaptive confidence interval for $X_{n+1}\theta$ is equivalent to having a fully adaptive confidence set for $\|\theta\|_2$. So, loosely speaking, doing inference on $L^T\theta$ for 'fully-loaded' isotropic linear functionals $L$ becomes as difficult as doing inference on $\theta$ in $\ell_2$. As a matter of fact, for

$$\mathcal{L} = diag(l_1, \ldots, l_j, \ldots, l_p), l_j \in \mathbb{R}$$

it can be shown that an adaptive confidence set for random linear functional $X^{n+1}\mathcal{L}\theta$ is equivalent to an adaptive confidence set under the loss function $\|\vartheta\|_{\mathcal{L}}$ defined by

$$\|\vartheta\|_{\mathcal{L}}^2 := \sum_{j=1}^{p}(l_j\vartheta_j)^2.$$

The proof is straightforward so we do not discuss it here. The upper bound in Theorem 2.2 suggests that for the class of random linear functionals with sub-Gaussian entries scaling to $(1^{-1/4}, 2^{-1/4}, \ldots, p^{-1/4})$, honest adaptive confidence intervals can be constructed. The lower bound displays, conceptually, the information-theoretic barrier for adaptation in regions where the sparsity of the parameter is significantly smaller than the cardinality of $L$, which is described as the 'dense' linear functional case in [10].

In practical contexts where the task is to identify 'anomalies' in the outcome, adaptive confidence sets for the linear functional can be immensely useful. For example in the field of quantitative finance where $Y$ is historical stock prices and $X$ describes the corresponding market conditions, $\theta$ is typically referred to as the 'factors' that drive stock price movements, which are high-dimensional. In such applications it is generally assumed that $\theta$ is sparse to some extent but specifying the range of sparsity may not be feasible. Given a new observation $(Y_{n+1}, X_{n+1})$ one may be interested in flagging unusual $Y_{n+1}$ for potential insider trading investigations, where $X_{n+1}^T \in \mathbb{R}^p$ is not necessarily sparse. Then $X_{n+1}\mathcal{L}\theta$ with $\mathcal{L} := diag(1^{-1/4}, 2^{-1/4}, \ldots, p^{-1/4})$ is the most 'informative' linear functional for which adaptive uncertainty quantification is guaranteed without further assumptions on $X_{n+1}$.

# 2 Main Results

## 2.1 Conditions and notations

We analyze the statistical model (1), assuming $k, p, n \to \infty$ and $p \simeq n$. We employ some standard conditions for inference in sparse linear models:

A. The design and noise: all $X_{ij}$ and $\varepsilon_i$ are i.i.d $N(0,1)$
B. The true parameter $\theta$:
   B1. Sparsity: $\|\theta\|_0 = k$, $k \leq k_{max}$ for some $k_{max} = o(n/\log p)$
   B2. The 'signal strength' is bounded: $\|\theta\|_2^2 := \sum_{j=1}^p \theta_j^2 \leq b^2$ for some known $0 < b < \infty$

Technical procedures for more general conditions are extensively studied in previous works mentioned in the introduction, especially in [7] and Chapters 2,6 and 7 in [2].

A *test* is any measurable function of $(Y, X)$ that takes values in $\{0, 1\}$.

We denote by $\mathbb{P}_\theta$ the law of $(Y, X)$ from model (1), and $\mathbb{E}_\theta$ the corresponding expectation. The expectation of $X$ is denoted by $E^X$ and $E_\theta$ is the expectation conditional on $X$. In this context when we say 'with high probability' we mean that under the law $\mathbb{P}_\theta$ the probability of a certain event becomes greater than $1-\epsilon$, for any $\epsilon > 0$, as $n, p \to \infty$.

For $\alpha \geq 0, \vartheta \in \mathbb{R}^p$, define norms

$$\|\vartheta\|_{\ell_2^{-\alpha}}^2 := \sum_{j=1}^p (j^{-\alpha}\vartheta_j)^2 \tag{6}$$

and when $\alpha = 0$ this becomes $\|\vartheta\|_2$.

For $\vartheta, \vartheta' \in \mathbb{R}^p$, let $(j)$ denote the reordering of coordinates such that $|\vartheta_{(1)}| \geq |\vartheta_{(2)}| \geq \cdots \geq |\vartheta_{(p)}|$, and define

$$d_{(-\alpha)}^2(\vartheta', \vartheta) := \sum_{j=1}^p j^{-2\alpha}(\vartheta_{(j)} - \vartheta'_{(j)})^2 \tag{7}$$

which we refer to as the $d_{(-\alpha)}(\cdot, \vartheta)$ -distance. The 'true' $d_{(-\alpha)}$ -distance refers to the distance based on the ordering of true values of $|\theta_j|$, namely $d_{(-\alpha)}(\cdot, \theta)$.

For $\vartheta \in \mathbb{R}^p$, denote by $B_0(k) := \{\vartheta : \|\vartheta\|_0 \leq k\}$ the '$l_0$'-ball that contains all $k$-sparse $\vartheta$. For $b > 0$, $q > 0$ let $B_0^q(k, b) := B_0(k) \cap \{\vartheta : \|\vartheta\|_q^q \leq b^q\}$, where $\|\vartheta\|_q$ is the usual $\ell_q$-norm.

For infinite sequences $y_n, x_n$, we use $x_n \lesssim y_n$ or $x_n = O(y_n)$ to imply that there exists C ¿ 0 such that $x_n \leq C y_n$ for all $n$ large. We write $x_n \simeq y_n$ if $x_n \lesssim y_n$ and $y_n \lesssim x_n$. We write $x_n = o(y_n)$ if $x_n/y_n \to 0$. $O_P(\cdot)$ and $o_P(\cdot)$ have their usual meanings of stochastic boundedness and convergence under the measure $P$.

For a semi-metric $\mathfrak{D}$, the $\mathfrak{D}$-diameter of a set $S$ is defined by $|S|_{\mathfrak{D}} := \sup\{\mathfrak{D}(x, y) : x, y \in S\}$, and the $\mathfrak{D}$-distance between sets $S_1, S_2$ is given by $\mathfrak{D}(S_1, S_2) := \inf\{\mathfrak{D}(x, y) : x \in S_1, y \in S_2\}$.

## 2.2 Adaptive estimation

Adaptive estimation of $\theta$ in model (1) under the $\ell_2$ and $l_\infty$ losses is extensively studied and known to be possible via numerous methods.

Let $\theta$ denote the true parameter and $\tilde{\theta}$ be any estimator that satisfies

$$\sup_{\theta \in B_0(k_{max})} \|\tilde{\theta} - \theta\|_\infty \leq \overline{C}\sqrt{\frac{\log p}{n}} \tag{8}$$

with high $\mathbb{P}_\theta$ probability for some constant $\overline{C}$. Under our assumptions the basic estimator $\tilde{\theta} = \frac{X^T Y}{n}$ possesses property (8), as $\frac{X^T Y}{n} - \theta = \frac{X^T \varepsilon}{n}$ and the supremum norm is taken over $p$ independent sub-exponential random variables. Under more general hypotheses for $(Y, X)$ such estimators also exist, such as the family of de-sparsified estimators discussed in [7], [8] and [9].

For $\alpha > \frac{1}{2}$, adaptive estimation is possible since, given (8),

$$\|\tilde{\theta} - \theta\|_{\ell_2^{-\alpha}}^2 = \sum_{j=1}^p |\tilde{\theta}_j - \theta_j|^2 \cdot j^{-2\alpha} \leq \|\tilde{\theta} - \theta\|_\infty^2 \cdot \sum_{j=1}^p j^{-2\alpha} \lesssim \frac{\log p}{n}$$

uniformly for all $\theta \in B_0(k_{max})$. So for $\alpha > \frac{1}{2}$ estimation in model 1 becomes a 'near-parametric' problem and the estimation rate does not depend on the underlying sparsity.

For $0 \le \alpha \le \frac{1}{2}$, thresholding any estimator that possesses property (8) suffices to achieve adaptive estimation:

**Theorem 2.1.** *Under model* (1) *and Condition A, let $\tilde{\theta}$ be any estimator that satisfies* (8)*. Let $\hat{\theta}_j = \tilde{\theta}_j \cdot \phi_j, 1 \le j \le p$ where*

$$\phi_j := \mathbf{1}\left\{\Big|\Big(\frac{X^T Y}{n}\Big)_j\Big| > C\sqrt{\frac{\log p}{n}}\right\} \tag{9}$$

*for some $0 < C < 1$.*

*Then $\forall \theta \in B_0(k)$, $k < k_{max}$ and any $0 \le \alpha \le \frac{1}{2}$, with high probability,*

$$\|\hat{\theta} - \theta\|^2_{\ell_2^{-\alpha}} \lesssim \frac{\log p}{n} k^{(1-2\alpha)}$$

A minimax lower bound for the estimation rate in $\|\cdot\|_{\ell_2^{-\alpha}}$ can be established in the same fashion as in $\ell_2$, through the well known procedure developed in [15] that reduces minimax estimation problems to the problem of multi-testing (see also Theorem 6.3.2 in [6]). The $\ell_2$ distance between hypotheses that lead to the lower bound is of the order $\sqrt{n^{-1} k \log p}$, which becomes $\sqrt{n^{-1} k^{(1-2\alpha)} \log p}$ in $\ell_2^{-\alpha}$. We refrain from going into details here.

### 2.3 Honest adaptive confidence sets

Having established the existence of adaptive estimators and the minimax estimation rate, we move onto showing our main result, namely that fully adaptive confidence sets exist in $\|\cdot\|_{\ell_2^{-\alpha}}$ when $\frac{1}{4} \le \alpha \le \frac{1}{2}$. We further show that the non-existence result for $\alpha = 0$ discussed in [3] extends to $0 < \alpha < 1/4$, at least if one is interested in adapting to the whole range of $1 \le k \le k_{max}$.

**Theorem 2.2.** *In model* (1)*, assume for all $1 \le i \le 2n$, $1 \le j \le p$, each $X_{ij}$ and $\varepsilon_i$ follows an independent $N(0,1)$ distribution and that the true parameter $\theta$ satisfies Condition B. Let $b < \infty$ be given, $n, k_{max} \to \infty$, $p \simeq n$, while $k_{max} = o(\frac{n}{\log p})$.*

*The upper bound: For $\frac{1}{4} \le \alpha \le \frac{1}{2}$, let $\hat{\theta}$ be as described in Theorem 2.1 based on samples $n+1 \le i \le 2n$ and $U_n^{-\alpha}(\hat{\theta})$ be defined as in* (5) *based on the first half of the sample. For any $\beta > 0$ and $n \in \mathbb{N}$ there exists $\mu_\beta$ such that the set*

$$C_n^{-\alpha} := \left\{\theta : \|\theta - \hat{\theta}\|^2_{\ell_2^{-\alpha}} \le U_n^{-\alpha}(\hat{\theta}) + \mu_\beta \frac{\log p}{n}\right\}$$

*satisfies*

$$\inf_{\theta \in B_0^2(k_{max}, b)} \mathbb{P}_\theta(\theta \in C_n^{-\alpha}) \ge 1 - \beta \tag{10}$$

*while for every $\beta' > 0$ there exists a universal constant $D$ s.t.*

$$\limsup_{n\to\infty} \sup_{\theta \in B_0^2(k,b)} \mathbb{P}_\theta\Big(|C_n^{-\alpha}|^2_{\ell_2^{-\alpha}} > D\frac{k^{(1-2\alpha)} \log p}{n}\Big) \le \beta' \tag{11}$$

*for any $1 \le k \le k_{max}$.*

*The lower bound: For $0 \le \alpha \le \frac{1}{2}$, let $0 < \beta, \beta' < 1$ be given. If for any $k_{max} \to \infty, k_{max} = o(\frac{n}{\log p})$ a confidence set $C_n^{-\alpha}$ satisfies*

$$\inf_{\theta \in B_0(k_{max})} \mathbb{P}_\theta(\theta \in C_n^{-\alpha}) \ge 1 - \beta \tag{12}$$

*for all $n$ large, then $C_n^{-\alpha}$ cannot also satisfy that*

$$\limsup_{n\to\infty} \sup_{\theta \in B_0(k)} \mathbb{P}_\theta\Big(|C_n^{-\alpha}|_{\ell_2^{-\alpha}} > r_n\Big) \le \beta' \tag{13}$$

*for all $1 \le k \le k_{max}$ at any rate*

$$r_n = o\left(\min\left\{\frac{p^{(\frac{1}{4}-\alpha)}}{\sqrt{n}}, \sqrt{\frac{k_{max}^{1-2\alpha}\log p}{n}}\right\}\right)$$

*In particular, when $0 \le \alpha < \frac{1}{4}$, there exists $k_{max} = o(\frac{\log p}{n})$ such that no confidence set can be adaptive over the range $1 \le k \le k_{max}$.*

The lower bound suggests that when $0 \le \alpha < \frac{1}{4}$ no honest confidence set can satisfy (13) with

$$r_n = \sqrt{\frac{k^{1-2\alpha}\log p}{n}},$$

which extends previous findings from $\ell_2$ in [3].

Perhaps surprisingly, in $\ell_2^{-\alpha}$ honest adaptive confidence sets can be constructed when $\alpha \ge \frac{1}{4}$. In practical applications it seems relevant if one could establish similar results for the *ordered* $\ell_2^{-\alpha}$-distance, where higher weighting is assigned to coordinates of $\theta$ with larger absolute value. The difficulty of constructing such adaptive confidence sets lies in the fact that the true parameter $\theta$, and thus the ordering, are not directly observable. Using techniques similar to those used in proving Theorem 2.2, one can show that ordering weightings by $|\hat{\theta}_j|$ in the $U$-statistic suffices.

Denote by $\Pi_p := (\pi_1, \pi_2, \dots, \pi_p)$ some permutation of the sequence $(1, 2, \dots, p)$. Let
$\bar{\Pi}_l = \{\pi_1, \pi_2, \dots, \pi_{l-1}\}$ while $\bar{\Pi}_1 = \emptyset$. Call the set $\bar{J}_l = \{1, 2, \dots, p\} \setminus \bar{\Pi}_l$.

For $\theta \in \mathbb{R}^p$ and an estimator $\hat{\theta}$ that is defined in Theorem 2.1, set the permutation $\hat{\Pi} \equiv (\hat{\pi}_1, \hat{\pi}_2, \dots, \hat{\pi}_p)$ as:

$$\begin{aligned}
\hat{\pi}_1 &= \arg\max_{j \in \bar{J}_1} |\hat{\theta}_j| \\
\hat{\pi}_2 &= \arg\max_{j \in \bar{J}_2} |\hat{\theta}_j| \\
&\vdots \\
\hat{\pi}_p &= \arg\max_{j \in \bar{J}_p} |\hat{\theta}_j|
\end{aligned}$$

Then for some estimator $\hat{\theta}$, in equivalence to definition (7), define for any $\vartheta \in \mathbb{R}^p$

$$d_{(-\alpha)}^2(\vartheta, \hat{\theta}) := \sum_{j=1}^{p} j^{-2\alpha}(\vartheta_{\hat{\pi}_j} - \hat{\theta}_{\hat{\pi}_j})^2$$

It is of course a random 'distance' as $\hat{\Pi}$ is random. But regardless of the permutation $\hat{\Pi}$, we measure the size of any confidence set by its $d_{(-\alpha)}$-diameter, i.e. under the ordering of the true values of $|\theta_j|$:

*Corollary* 2.3 (Adaptive inference in the ordered distance). Let $\hat{\theta}$ be defined as in Theorem 2.1 and based on an independent sample of size $n$ generated from $\theta$. Given $\hat{\theta}$ and some $(Y, X)$ under model (1), let the permutation $\hat{\Pi}$ be as above. For $\vartheta \in \mathbb{R}^p$, define

$$U_n^{(-\alpha)}(\vartheta) := \frac{2}{n(n-1)} \sum_{j=1}^{p} \sum_{l=2}^{n} \sum_{m=1}^{l-1} (Y_l X_{l\hat{\pi}_j} - \vartheta_{\hat{\pi}_j})(Y_m X_{m\hat{\pi}_j} - \vartheta_{\hat{\pi}_j}) \cdot j^{-2\alpha}$$

Then for $\frac{1}{4} \le \alpha \le \frac{1}{2}$ and any $\beta > 0$, there exists $\mu'_\beta$ such that the set

$$C_n^{(-\alpha)} := \left\{\theta : d_{(-\alpha)}^2(\theta, \hat{\theta}) \le U_n^{(-\alpha)}(\hat{\theta}) + \mu'_\beta \frac{\log p}{n}\right\} \tag{14}$$

satisfies that

$$\inf_{\theta \in B_0^2(k_{max}, b)} \mathbb{P}_\theta(\theta \in C_n^{(-\alpha)}) \ge 1 - \beta$$

while for every $\beta' > 0$ there exists a universal constant $D$ s.t.

$$\limsup_{n\to\infty} \sup_{\theta\in B_0^2(k,b)} \mathbb{P}_\theta\Big(|C_n^{(-\alpha)}|^2_{d_{(-\alpha)}} > D\frac{k^{(1-2\alpha)}\log p}{n}\Big) \leq \beta'$$

for any $1 \leq k \leq k_{max}$

Note that, recalling definition (7), the confidence set $C_n^{(-\alpha)}$ is defined with the $d_{(-\alpha)}(\cdot,\hat{\theta})$-distance, but it has the desired diameter in the true ordered distance $d_{(-\alpha)}(\cdot,\theta)$ nonetheless.

# 3 Proofs for main results

## 3.1 Proof of Theorem 2.1

*Proof.* In model (1) and under our conditions, as $X_{ij}, \varepsilon_j \overset{\text{i.i.d.}}{\sim} N(0,1)$:

$$\mathbb{E}_\theta\left[\Big(\frac{X^TY}{n}\Big)_j\right] = \mathbb{E}_\theta\left[\Big(\frac{X^TX\theta}{n}\Big)_j + \Big(\frac{X^T\varepsilon}{n}\Big)_j\right] = \theta_j$$

When $\theta_j \geq \sqrt{\frac{\log p}{n}}$, for any $0 < C < 1$,

$$\begin{aligned}
\mathbb{P}_\theta\left(\Big|\Big(\frac{X^TY}{n}\Big)_j\Big| > C\sqrt{\frac{\log p}{n}}\right) &\geq \mathbb{P}_\theta\left(\Big(\frac{X^TY}{n}\Big)_j > C\sqrt{\frac{\log p}{n}}\right) \\
&= \mathbb{P}_\theta\left(\Big(\frac{X^T\varepsilon}{n}\Big)_j + \Big(\frac{X^TX\theta}{n}\Big)_j - \mathbb{E}_\theta\Big[\Big(\frac{X^TX\theta}{n}\Big)_j\Big] > (C-1)\sqrt{\frac{\log p}{n}}\right) \\
&\geq 1 - \mathbb{P}_\theta(\mathbb{A}) - \mathbb{P}_\theta(\mathbb{B})
\end{aligned}$$

where

$$\mathbb{P}_\theta(\mathbb{A}) := \mathbb{P}_\theta\left(\Big(\frac{X^TX\theta}{n}\Big)_j - \mathbb{E}_\theta\Big[\Big(\frac{X^TX\theta}{n}\Big)_j\Big] < \frac{C-1}{2}\sqrt{\frac{\log p}{n}}\right)$$

$$\mathbb{P}_\theta(\mathbb{B}) := \mathbb{P}_\theta\left(\Big(\frac{X^T\varepsilon}{n}\Big)_j < \frac{C-1}{2}\sqrt{\frac{\log p}{n}}\right)$$

Using Lemma A.2 and Lemma A.4 in Appendix A, both probabilities approach 0 as $p \to \infty$. Combining the two probabilities there is, for suitable choices of $C, C_1$,

$$\mathbb{P}_\theta\left(\phi_j = 0 \middle| \theta_j \geq \sqrt{\frac{\log p}{n}}\right) \leq C_1\exp(-\log p)$$

Similarly, using $\mathbb{P}_\theta\left(\Big|\Big(\frac{X^TY}{n}\Big)_j\Big| > C\sqrt{\frac{\log p}{n}}\right) \geq \mathbb{P}_\theta\left(\Big(\frac{X^TY}{n}\Big)_j < -C\sqrt{\frac{\log p}{n}}\right)$ we can show that

$$\mathbb{P}_\theta\left(\phi_j = 0 \middle| \theta_j \leq -\sqrt{\frac{\log p}{n}}\right) \leq C_1\exp(-\log p)$$

So, by assumptions on $p, k$ and when $n \to \infty$,

$$\sum_{j:\phi_j=0} \mathbf{1}\left(|\theta_j| \geq \sqrt{\frac{\log p}{n}}\right) = 0 \tag{15}$$

with high probability.

When $\theta_j = 0$, again using Lemma A.2 to bound $X^T\varepsilon$ we have:

$$\mathbb{P}_\theta\Big(\phi_j = 1\Big|\theta_j = 0\Big) = \mathbb{P}_\theta\left(\Big(\frac{X^T\varepsilon}{n}\Big)_j < (C-1)\sqrt{\frac{\log p}{n}}\right) \le C_2\exp(-C_3\log p) \tag{16}$$

for any $C_3 > 0$ and a suitable choice of $C_2$.

Now we want to control $\sum_{i=1}^{p}\phi_j$. For $\theta$ at most $k$-sparse, take any $C_3 > 2$:

$$\begin{aligned}
\mathbb{P}_\theta\Big(\sum_{j=1}^{p}\phi_j > 2k\Big) &\le \mathbb{P}_\theta\Big(\sum_{j:\theta_j=0}\phi_j > k\Big) \\
&\le \sum_{l=k}^{p}\binom{p}{l}\Big(\mathbb{P}_\theta\big(\phi_j = 1|\theta_j = 0\big)\Big)^l \\
&< \sum_{l=k}^{p}(C_2\cdot p^{1-C_3})^l < p\cdot(C_2\cdot p^{1-C_3})^k \\
&< C_2\cdot p^{2-C_3}
\end{aligned} \tag{17}$$

which can be arbitrarily small given choice of $C_3$.

We now show that the estimation error is controlled in $\ell_2^{-\alpha}$ regardless of how the coordinates are ordered:

**Lemma 3.1.** *For $0 \le \alpha \le \frac{1}{2}$ and any permutation $\Pi_p := (\pi_1, \dots, \pi_j, \dots, \pi_p)$ of the sequence $(1, 2, \dots, p)$, for some constant independent of $\Pi_p$*

$$\sum_{j=1}^{p}(\tilde{\theta}_{\pi_j}\cdot\phi_{\pi_j} - \theta_{\pi_j})^2\cdot j^{-2\alpha} \lesssim \frac{\log p}{n}k^{(1-2\alpha)}$$

*with high probability.*

The above lemma can be shown using (8), (15), (17) and the fact that $\theta \in B_0(k)$:

$$\begin{aligned}
\sum_{j=1}^{p}(\tilde{\theta}_{\pi_j}\cdot\phi_{\pi_j} - \theta_{\pi_j})^2\cdot j^{-2\alpha} &= \sum_{j:\phi_{\pi_j}=1}(\tilde{\theta}_{\pi_j} - \theta_{\pi_j})^2\cdot j^{-2\alpha} + \sum_{j:\phi_{\pi_j}=0}\theta_{\pi_j}^2\cdot j^{-2\alpha} \\
&\le \|\tilde{\theta} - \theta\|_\infty^2\sum_{j:\phi_{\pi_j}=1}j^{-2\alpha} + \sum_{j:\theta_{\pi_j}^2\le\frac{\log p}{n}}\theta_{\pi_j}^2\cdot j^{-2\alpha} \\
&\quad + \sum_{j:\phi_{\pi_j}=0}\theta_{\pi_j}^2\mathbf{1}\left(|\theta_{\pi_j}| \ge \sqrt{\frac{\log p}{n}}\right)j^{-2\alpha} \\
&\lesssim \frac{\log p}{n}\Big(\sum_{j:\phi_{\pi_j}=1}j^{-2\alpha} + \sum_{j:0<\theta_{\pi_j}^2\le\frac{\log p}{n}}j^{-2\alpha}\Big) \lesssim \frac{\log p}{n}k^{(1-2\alpha)}
\end{aligned}$$

with high probability.

Taking $\Pi = (1, 2, \dots, p)$ leads to Theorem 2.1.

□

### 3.2 Proof of Theorem 2.2 - The upper bound: an adaptive confidence set based on the U-statistic

*Proof.* Without loss of generality we can assume the sample size is $2n$. Using the first half of the sample and an estimator that satisfies (8), we can construct $\hat{\theta} = \hat{\theta}_n$ s.t. for any $k = o(n/\log p)$

$$\|\hat{\theta} - \theta\|^2_{\ell_2^{-\alpha}} \lesssim \frac{\log p}{n} k^{(1-2\alpha)} \tag{18}$$

with high $\mathbb{P}_\theta$-probability, as shown in Theorem 1.

We denote the second half of the sample $(Y, X)$ where $Y \in \mathbb{R}^n$ and $X \in \mathbb{R}^{n \times p}$. As specified by our conditions, $X_{ij}, \varepsilon_j \overset{\text{i.i.d.}}{\sim} N(0, 1)$. Below we take expectations w.r.t. $(Y, X)$ only.

For $\vartheta \in \mathbb{R}^p$, define the $U$-statistic

$$U_n^{-\alpha}(\vartheta) := \frac{2}{n(n-1)} \sum_{j=1}^{p} \sum_{l=2}^{n} \sum_{m=1}^{l-1} (Y_l X_{lj} - \vartheta_j)(Y_m X_{mj} - \vartheta_j) \cdot j^{-2\alpha} \tag{19}$$

We show below that, with the suitbale choice of $\mu_\beta$ as discussed, the confidence set

$$C_n^{-\alpha} := \left\{ \theta : \|\theta - \hat{\theta}\|^2_{\ell_2^{-\alpha}} \leq U_n^{-\alpha}(\hat{\theta}) + \mu_\beta \frac{\log p}{n} \right\}$$

achieves the desired diameter and coverage without requiring knowledge of the true sparsity.

As $\mathbb{E}_\theta[Y_i X_{ij}] = \theta_j$ for any $i$,

$$\mathbb{E}_\theta[U_n^{-\alpha}(\hat{\theta})] = \|\theta - \hat{\theta}\|^2_{\ell_2^{-\alpha}}$$

By the Hoeffding Decomposition [16] we have

$$U_n^{-\alpha}(\hat{\theta}) = U_n^{-\alpha}(\theta) + 2L_n^{-\alpha}(\hat{\theta}) + \|\hat{\theta} - \theta\|^2_{\ell_2^{-\alpha}} \tag{20}$$

where

$$L_n^{-\alpha}(\hat{\theta}) := \frac{1}{n} \sum_{i=1}^{n} \sum_{j=1}^{p} (\theta_j - Y_i X_{ij})(\theta_j - \hat{\theta}_j) \cdot j^{-2\alpha}$$

By (18) and Lemma 4.1, for $1/4 \leq \alpha < 1/2$ and $\theta \in B_0^2(k_{max}, b)$,

$$\begin{aligned} |C_n^{-\alpha}|_{\ell_2^{-\alpha}} &= 2\sqrt{U_n^{-\alpha}(\theta) + 2L_n^{-\alpha}(\hat{\theta}) + \|\hat{\theta} - \theta\|^2_{\ell_2^{-\alpha}} + \mu_\beta \frac{\log p}{n}} \\ &\lesssim \sqrt{\frac{k^{(1-2\alpha)} \log p}{n}} \end{aligned}$$

with high probability. This gives the desired diameter in the sense of (11).

For any $\theta \in B_0^2(k_{max}, b)$ the coverage probability

$$\begin{aligned} \mathbb{P}_\theta\{\theta \notin C_n^{-\alpha}\} &= \mathbb{P}_\theta\left\{ \|\theta - \hat{\theta}\|^2_{\ell_2^{-\alpha}} \geq U_n^{-\alpha}(\hat{\theta}) + \mu_\beta \frac{\log p}{n} \right\} \\ &= \mathbb{P}_\theta\left\{ U_n^{-\alpha}(\hat{\theta}) - \|\theta - \hat{\theta}\|^2_{\ell_2^{-\alpha}} \leq -\mu_\beta \frac{\log p}{n} \right\} \\ &= \mathbb{P}_\theta\left\{ U_n^{-\alpha}(\theta) + 2L_n^{-\alpha}(\hat{\theta}) \leq -\mu_\beta \frac{\log p}{n} \right\} \\ &\leq \mathbb{P}_\theta\left\{ |L_n^{-\alpha}(\hat{\theta})| \geq \frac{\mu_\beta}{4} \frac{\log p}{n} \right\} + \mathbb{P}_\theta\left\{ |U_n^{-\alpha}(\theta)| \geq \frac{\mu_\beta}{2} \frac{\log p}{n} \right\} \end{aligned}$$

By Lemma 4.1 and Chebyshev's inequality, both terms in the last inequality can be made small enough by a suitable choice of $\mu_\beta$. For the first term:

$$\mathbb{P}_\theta\left\{|L_n^{-\alpha}(\hat{\theta})| \geq \frac{\mu_\beta}{4}\frac{\log p}{n}\right\} \leq \frac{16}{\mu_\beta^2}\frac{n^2}{(\log p)^2}Var(L_n^{-\alpha}(\hat{\theta})) \leq \frac{16C_b}{\mu_\beta^2}.$$

And similarly for the second term. So we have 10 and the upper bound. □

### 3.3 Proof of Theorem 2.2 - The lower bound

*Proof.* Via general decision-theoretic techniques one can show that the existence of an honest adaptive confidence set implies that one is able to consistently test between different sub-models of the parameter space. More specifically, for the following composite testing problem:

$$H_0 : \theta \in B_0(k') \quad vs \quad H_1 : \theta \in \widetilde{B}_0^\alpha(k,\rho) \tag{21}$$

where $k \leq k_{max}$, $k' = o(k^{1-2\alpha})$ as $n, k \to \infty$ and

$$\widetilde{B}_0^\alpha(k,\rho) := \{\theta : \|\theta\|_0 \leq k, \inf_{\vartheta \in B_0(k')} \|\theta - \vartheta\|_{\ell_2^{-\alpha}} \geq \rho\}.$$

If there exists some $C_n^{-\alpha}$ satisfying (12), and for some infinite sequence $\rho_n$

$$\sup_{\theta \in B_0(k')} \mathbb{E}_\theta[|C_n^{-\alpha}|_{\ell_2^{-\alpha}}] \lesssim \rho_n,$$

then the test

$$\psi_n = \mathbf{1}\{C_n^{-\alpha} \cap \tilde{B}_0^\alpha(k,\rho_n) \neq \emptyset\}$$

is a consistent test for (21) at rate $\rho_n$. See Proposition 8.3.6 in [6] for a detailed account of this equivalence.

We now show that, when $p \simeq n$ and $0 \leq \alpha \leq \frac{1}{2}$,

$$\rho_{n,p,k} \simeq \min\left\{\sqrt{n^{-1}k_{max}^{(1-2\alpha)}\log p}, p^{(\frac{1}{4}-\alpha)}n^{-\frac{1}{2}}\right\}$$

is a lower bound for the minimax testing rate for (21), i.e. any $\rho^* \equiv \rho_{n,p}^* = o(\rho_{n,p,k})$ leads to consistent testing being impossible. When $0 \leq \alpha < \frac{1}{4}$, the estimation rate in lower-dimensional sub-models becomes very fast while sufficient separation $\rho \gtrsim p^{(\frac{1}{4}-\alpha)}n^{-\frac{1}{2}}$ is still required for testing to be possible, resulting in the lower bound.

The proof uses the following key technical lemma:

**Lemma 3.2.** *For testing problem* (21)*, suppose there exists a sequence of probability measures $\pi = \pi_{n,p}$ on $\mathbb{R}^p$ such that $\pi(H_1) \to 1$ as $n, p \to \infty$. For $\vartheta, \vartheta' \sim \pi_{n,p}$ independently, if $\|\vartheta\|_2 = o_p(n^{-\frac{1}{4}})$ under $\pi$ and*

$$\mathbb{E}_{\pi^2}\left[\exp(n\vartheta^T\vartheta')\right] \leq 1,$$

*then we necessarily have*

$$\liminf_{n,p,k} \inf_\psi \Big(\sup_{\vartheta \in H_0} \mathbb{E}_\vartheta[\psi] + \sup_{\vartheta \in H_1} \mathbb{E}_\vartheta[1-\psi]\Big) \geq 1 \tag{22}$$

*where the infimium is taken over all tests $\psi$ for* (21)*.*

The inequality (22) essentially suggests that no effective test exists for testing problem (21). The techniques for proving the above lemma were detailed in [4]. Also refer to [17], [18] and [19] for proofs in settings similar to our problem. For completeness an independent proof is laid out the supplement material.

Now we construct a $\pi$ that satisfies conditions specified in Lemma 3.2:

With $n, p, k \to \infty$, let

$$h = \frac{ck^{(1-2\alpha)}\log p}{\sum_{m=1}^p m^{-2\alpha}}, t = \frac{\rho^*}{c\sqrt{k^{(1-2\alpha)}\log p}}$$

for some $0 < c < 1$. Consider a random vector $\vartheta \in \mathbb{R}^p$ such that each component $\vartheta_m$ of $\vartheta$ independently takes value $\vartheta_m = t\epsilon_m$, where $\epsilon_m \in \{-1, 0, 1\}$ with probability

$$\mathbb{P}(\epsilon_m = 0) = 1 - h, \text{ and } \mathbb{P}(\epsilon_m = 1) = \mathbb{P}(\epsilon_m = -1) = \frac{h}{2}. \tag{23}$$

This defines a probability measure on $\vartheta_m$ which we call $\pi_m$. Take $\pi$ as the product measure $\pi = \prod_{m=1}^{p} \pi_m$. Notice that when $\alpha = 0$ this becomes the standard $\ell_2$ scenario and the proof below leads to the lower bound presented in [3].

First we check that for some $k = o(p)$ and $H_1$ as defined in (21), $\pi(H_1) \to 1$.

The random variable $\Lambda := \sum_{m=1}^{p} |\epsilon_m|$ follows a binomial distribution where $\mathbb{E}[\Lambda] = ph$ and $Var(\Lambda) = ph(1-h)$. Using Chebyshev's Inequality:

$$\begin{aligned}
\pi(\vartheta \notin B_0(k)) &= \mathbb{P}\Big(\sum_{m=1}^{p} |\epsilon_m| > k\Big) \\
&= \mathbb{P}\Big(\sum_{m=1}^{p} (|\epsilon_m| - \mathbb{E}[|\epsilon_m|]) > (k - ph)\Big) \\
&\leq \frac{ph(1-h)}{(k-ph)^2}
\end{aligned} \tag{24}$$

As $h = o(1)$, there is $k = o(ph)$ and $\frac{ph(1-h)}{(k-ph)^2} \simeq \frac{1}{ph} = o(1)$. So the above probability $\pi(\vartheta \notin B_0(k)) \to 0$.

If $\vartheta \in B_0(k)$, then the infimium

$$\inf_{\vartheta' \in B_0(k')} \|\vartheta - \vartheta'\|_{\ell_2^{-\alpha}} = \|\vartheta - \vartheta^*\|^2_{\ell_2^{-\alpha}}$$

where $\vartheta^*$ is attained by setting the $k'$ coordinates in $\vartheta$ with the largest $|\vartheta_m| \cdot m^{-\alpha}$ values to 0.

Let $W_m = m^{-\alpha}\varepsilon_m^2$ and $W_{(m)}$ are $W_m$ reordered so $W_{(1)} < W_{(2)} < \cdots < W_{(p)}$.

Note also for $k' = o(k^{1-2\alpha})$, we have

$$t^2 k' = \frac{(\rho^*)^2}{c^2 k^{1-2\alpha} \log p} \times k' = o((\rho^*)^2)$$

Then $\forall \delta \in (0,1)$ and $k'$ small enough:

$$\begin{aligned}
\pi\Big(\inf_{\vartheta' \in B_0(k')} \|\vartheta - \vartheta'\|_{\ell_2^{-\alpha}} < (1-\delta)(\rho^*)^2\Big) &= \mathbb{P}(b^2 \sum_{m=1}^{p-k'} W_{(m)} < (1-\delta))(\rho^*)^2) \\
&= \mathbb{P}(t^2 \sum_{m=1}^{p-k'} W_{(m)} + b^2 k' < (1-\delta))(\rho^*)^2 + t^2 k') \\
&\leq \mathbb{P}(t^2 \sum_{m=1}^{p} W_m < (\rho^*)^2) \\
&= \pi(\|\theta\|^2_{\ell_2^{-\alpha}} < (\rho^*)^2)
\end{aligned}$$

Under $\pi$ we also have:

$$\mathbb{E}\Big[\|\vartheta\|^2_{\ell_2^{-\alpha}}\Big] = t^2 h \sum_{m=1}^{p} m^{-2\alpha} = \frac{1}{c}(\rho^*)^2 > (\rho^*)^2$$

So we deduce

$$\begin{aligned}\pi(\|\vartheta\|^2_{\ell_2^{-\alpha}} < (\rho^*)^2) &= \mathbb{P}\Big(\|\vartheta\|^2_{\ell_2^{-\alpha}} - \mathbb{E}[\|\vartheta\|^2_{\ell_2^{-\alpha}}] < (\rho^*)^2 - \mathbb{E}[\|\vartheta\|^2_{\ell_2^{-\alpha}}]\Big)\\ &= \mathbb{P}\Big(-t^2\sum_{m=1}^{p}(\epsilon_m^2 m^{-2\alpha} - \mathbb{E}[\epsilon_m^2 m^{-2\alpha}]) > \mathbb{E}[\|\vartheta\|^2_{\ell_2^{-\alpha}}] - (\rho^*)^2\Big)\\ &\leq \mathbb{P}\Big(|t^2\sum_{m=1}^{p}(\epsilon_m^2 m^{-2\alpha} - \mathbb{E}[\epsilon_m^2 m^{-2\alpha}])| > \mathbb{E}[\|\vartheta\|^2_{\ell_2^{-\alpha}}] - (\rho^*)^2\Big)\\ &= \mathbb{P}\Big(|t^2\sum_{m=1}^{p}(\epsilon_m^2 m^{-2\alpha} - \mathbb{E}[\epsilon_m^2 m^{-2\alpha}])| > (\frac{1}{c}-1)(\rho^*)^2\Big)\end{aligned}$$

The random variable $\Lambda' \coloneqq \sum_{m=1}^{p}\epsilon_m^2 m^{-2\alpha}$ follows a binomial distribution where $\mathbb{E}[\Lambda'] = h\sum_{m=1}^{p}m^{-2\alpha}$ and $Var(\Lambda') = h(1-h)\sum_{m=1}^{p}m^{-4\alpha}$. Using Chebyshev's Inequality, the last probability is bounded above by:

$$\begin{aligned}&\mathbb{P}(|\Lambda' - \mathbb{E}[\Lambda']| > (\frac{1}{c}-1)\frac{1}{t^2}(\rho^*)^2\frac{1}{\sqrt{h(1-h)\sum_{m=1}^{p}m^{-4\alpha}}}\sigma(\Lambda'))\\ &\leq \frac{c^2t^4h(1-h)\sum_{m=1}^{p}m^{-4\alpha}}{(1-c)^2(\rho^*)^4}\\ &= O(\frac{1}{k^{1-2\alpha}p^{2\alpha}}) \to 0\end{aligned}$$

So $\pi(H_1) = \pi_n\big(\vartheta \in B_0(k), \inf_{\vartheta' \in B_0(k')}\|\vartheta - \vartheta'\|_{\ell_2^{-\alpha}} > \rho^*\big) \to 1$.

Now we evaluate

$$\mathbb{E}_{\pi^2}\Big[\exp\Big(n\sum_{i=1}^{p}\vartheta_i \cdot \vartheta'_i\Big)\Big]$$

where $\vartheta$ and $\vartheta'$ are independent draws from $\pi$.

Let $\epsilon_m$, $\epsilon'_m$ be i.i.d with distribution specified in (23), and $V = \epsilon_m \cdot \epsilon'_m$. Then $\mathbb{E}[V] = 0$,

$$n\sum_{i=1}^{p}\vartheta_i \cdot \vartheta'_i = nt^2\sum_{m=1}^{p}\epsilon_m \cdot \epsilon'_m$$

and

$$\begin{aligned}Var(V) &= \mathbb{E}[V^2]\\ &= \mathbb{E}[V^2|V \neq 0]\cdot\mathbb{P}[V \neq 0] + \mathbb{E}[V^2|V = 0]\cdot\mathbb{P}[V = 0]\\ &= h^2\end{aligned}$$

As $|V| \leq 1$, we can use Theorem 3.1.5 from [6]:

$$\begin{aligned}\mathbb{E}_{\pi^2}\Big[\exp n\sum_{i=1}^{p}\vartheta_i \cdot \vartheta'_i\Big] &= \mathbb{E}_{\pi^2}\Big[\exp nt^2\sum_{m=1}^{p}\epsilon_m \cdot \epsilon'_m\Big]\\ &\leq \exp(ph^2(e^{nt^2} - 1 - nt^2))\end{aligned}$$

If $\rho^* = o\left(\sqrt{\frac{k^{(1-2\alpha)}\log p}{n}}\right)$ and $\rho^* = o\left(\frac{p^{(\frac{1}{4}-\alpha)}}{\sqrt{n}}\right)$ at the same time, we have $nt^2 = o(1)$ so

$$e^{nt^2} = 1 + nt^2 + \frac{(nt^2)^2}{2} + o((nt^2)^2)$$

and

$$\mathbb{E}_{\pi^2}\Big[\exp n\sum_{i=1}^{p}\vartheta_i\cdot\vartheta_i'\Big]\le\exp(ph^2n^2t^4(1+o(1)))$$

where

$$ph^2n^2t^4=O\Big(p\cdot\frac{k^{(2-4\alpha)}}{p^{(2-4\alpha)}}\cdot n^2\cdot\frac{(\rho^*)^4}{k^{(2-4\alpha)}}\Big)=o(1)$$

By Lemma 3.2, this leads to

$$\sup_{\vartheta\in H_0}\mathbb{E}_\vartheta[\psi]+\sup_{\vartheta\in H_1}\mathbb{E}_\vartheta[1-\psi]\ge 1$$

for all $\psi$. So testing (21) at any rate faster than $\rho_{n,p,k}$ is not possible.

We conclude the proof that, with the 'equivalence' between testing and honest adaptive confidence sets in mind, the $\ell_2^{-\alpha}$ diameter of any honest adaptive confidence set $C_n$ is at least $O_P(r_n)$ with

$$r_n=\min\left(\sqrt{\frac{k_{max}^{(1-2\alpha)}\log p}{n}},\frac{p^{(\frac{1}{4}-\alpha)}}{\sqrt{n}}\right)$$

as required.

□

### 3.4 Proof of Corollary 2.3

*Proof.* By Lemma 3.1, for any true parameter $\theta$, $\hat{\theta}$ is an adaptive estimator under $d_{(-\alpha)}(\cdot,\theta)$ distance.

First we check the diameter of $C_n^{(-\alpha)}$ in $d_{(-\alpha)}(\cdot,\theta)$:

Let

$$L_n^{(-\alpha)}(\hat{\theta}):=\frac{1}{n}\sum_{i=1}^{n}\sum_{j=1}^{p}(\theta_{\hat{\pi}_j}-Y_iX_{i\hat{\pi}_j})(\theta_{\hat{\pi}_j}-\hat{\theta}_{\hat{\pi}_j})\cdot j^{-2\alpha}.$$

By Hoeffding Decomposition (20):

$$\begin{aligned}C_n^{(-\alpha)}&=\left\{\theta:d^2_{(-\alpha)}(\theta,\hat{\theta})\le U_n^{(-\alpha)}(\hat{\theta})+\mu_\beta'\frac{\log p}{n}\right\}\\&=\left\{\theta:d^2_{(-\alpha)}(\theta,\hat{\theta})\le U_n^{(-\alpha)}(\theta)+2L_n^{(-\alpha)}(\hat{\theta})+d^2_{(-\alpha)}(\theta,\hat{\theta})+\mu_\beta'\frac{\log p}{n}\right\}\\&=\left\{\theta:d^2_{(-\alpha)}(\theta,\theta)\le U_n^{(-\alpha)}(\theta)+2L_n^{(-\alpha)}(\hat{\theta})+d^2_{(-\alpha)}(\theta,\theta)+\mu_\beta'\frac{\log p}{n}\right\}\end{aligned}$$

The variance of $U_n^{(-\alpha)}(\theta)$ and $L_n^{(-\alpha)}(\hat{\theta})$ can be bounded to the order of $\frac{(\log p)^2}{n^2}$ with same procedures as in Lemma 4.1. We also have $d^2_{(-\alpha)}(\theta,\theta)\lesssim\frac{k^{(1-2\alpha)}\log p}{n}$ from the fact that $\hat{\theta}$ is an adaptive estimator in $d_{(-\alpha)}(\cdot,\theta)$. So

$$|C_n^{(-\alpha)}|_{(-\alpha)}=O_P\left(\sqrt{\frac{k^{(1-2\alpha)}\log p}{n}}\right)$$

For coverage, again by Hoeffding Decomposition

$$\begin{aligned}
\mathbb{P}_\theta\{\theta \notin C_n^{(-\alpha)}\} &= \mathbb{P}_\theta\left\{d^2_{(-\alpha)}(\theta,\hat\theta) \geq U_n^{(-\alpha)}(\theta) + 2L_n^{(-\alpha)}(\hat\theta) + d^2_{(-\alpha)}(\theta,\hat\theta) + \mu'_\beta \frac{\log p}{n}\right\} \\
&= \mathbb{P}_\theta\left\{U_n^{(-\alpha)}(\theta) + 2L_n^{(-\hat\alpha)}(\hat\theta) \leq -\mu'_\beta \frac{\log p}{n}\right\} \\
&\leq \mathbb{P}_\theta\left\{|U_n^{(-\alpha)}(\theta)| + 2|L_n^{(-\hat\alpha)}(\hat\theta)| \geq \mu'_\beta \frac{\log p}{n}\right\} \\
&\leq \mathbb{P}_\theta\left\{|U_n^{(-\alpha)}(\theta)| \geq \frac{\mu'_\beta}{2} \frac{\log p}{n}\right\} + \mathbb{P}_\theta\left\{|L_n^{(-\hat\alpha)}(\hat\theta)| \geq \frac{\mu'_\beta}{4} \frac{\log p}{n}\right\},
\end{aligned}$$

and $\mu'_\beta$ can be chosen for any desired coverage probability $1-\beta$ through Chebyshev's inequality.

□

# 4 Technical proofs

*Proof of Lemma 3.2.* Note that proving for $H_0 = \{0\}$ is sufficient. Denote by $\mathbb{E}_0[\cdot]$ the expectation when $\theta = 0$.

Let $\pi = \pi_{n,p}$ be a sequence of probability measures on $\mathbb{R}^p$ such that $\pi(H_1) \to 1$ as $n, p \to \infty$.

We use $\mathbb{E}_{\theta\sim\pi}[\cdot]$ to denote the expectation under $\theta$ when $\theta$ is drawn from $\pi$. Also we shorthand $\mathbb{E}_{\pi^2}$ for $\mathbb{E}_{\theta\times\theta'\sim\pi\times\pi}$, where $\theta \times \theta' \sim \pi \times \pi$ means that $\theta$ and $\theta'$ are drawn independently from $\pi$. Then

$$\Big(\mathbb{E}_{\theta\sim\pi}[f(\theta)]\Big)^2 = \mathbb{E}_{\pi^2}[f(\theta)\cdot f(\theta')]. \tag{25}$$

Let $\pi|H_1$ denote $\pi$ restricted and re-normalized on $H_1$. That is, $\forall\theta \in H_1$

$$\mathbb{E}_{\theta\sim\pi|H_1}[\cdot] = \frac{1}{\pi(H_1)}\mathbb{E}_{\theta\sim\pi}[\cdot]$$

Denote

$$Z = \mathbb{E}_{\theta\sim\pi}\prod_{i\leq n}\frac{dP_i^{(\theta)}}{dP_i^{(0)}} = \int\prod_{i\leq n}\frac{dP_i^{(\theta)}}{dP_i^{(0)}}d\pi(\theta),$$

where $dP_i^{(\theta)}$ is the distribution of $Y|X$ when $Y$ is generated by some $\theta \in H_1$, and $dP_i^{(0)}$ is the distribution of $Y|X$ when $Y$ is generated by $\theta = 0$.

If $\pi(H_1) \to 1$, using Markov's inequality we have, for any $\eta \in (0,1)$ and any test $\psi$:

$$\begin{aligned}
\mathbb{E}_0[\psi] + \sup_{\theta\in H_1}\mathbb{E}_\theta[1-\psi] &\geq \mathbb{E}_0[\psi] + \mathbb{E}_{\theta\sim\pi|H_1}\mathbb{E}_\theta[1-\psi] \\
&\geq \mathbb{E}_0[\psi] + \mathbb{E}_{\theta\sim\pi}\mathbb{E}_\theta[1-\psi] - o(1) \\
&= E^X[E_0[\psi] + \mathbb{E}_{\theta\sim\pi}E_\theta[1-\psi]] - o(1) \\
&= E^X[E_0[\psi] + ZE_0[1-\psi]] - o(1) \\
&= \mathbb{E}_0[\mathbf{1}\{\psi = 1\} + \mathbf{1}\{\psi = 0\}Z] - o(1) \\
&\geq \mathbb{E}_0[(1-\eta)\mathbf{1}\{Z \geq 1-\eta\}] - o(1) \\
&\geq (1-\eta)\left[1 - \left(\frac{\mathbb{E}_0[Z-1]^2}{\eta^2}\right)^{\frac{1}{2}}\right] - o(1)
\end{aligned} \tag{26}$$

Given the model (1) we have:

$$\begin{aligned} Z &= \mathbb{E}_{\theta\sim\pi} \prod_{i\le n} \frac{dP_i^{(\theta)}}{dP_i^{(0)}} = \mathbb{E}_{\theta\sim\pi}\left[\prod_{i\le n} \frac{\exp(-\frac{1}{2}(Y_i-(X\theta)_i)^2)}{\exp(-\frac{1}{2}(Y_i)^2)}\right] \\ &= \mathbb{E}_{\theta\sim\pi}\Big[\prod_{i\le n} \exp\Big(\frac{1}{2}Y_i^2 - \frac{1}{2}(Y_i-(X\theta)_i)^2\Big)\Big] \\ &= \mathbb{E}_{\theta\sim\pi}\Big[\prod_{i\le n} \exp\Big(Y_i(X\theta)_i - \frac{1}{2}(X\theta)_i^2\Big)\Big] \end{aligned}$$

Easy to see that $\mathbb{E}_0[Z] = 1$, and

$$\mathbb{E}_0[(Z-1)^2] = \mathbb{E}_0[Z^2] - 1 \tag{27}$$

For n independent draws from $Y_i \overset{\text{i.i.d.}}{\sim} N((X\theta)_i, 1)$:

$$\begin{aligned} E_0[Z^2] &= \int_{\mathbb{R}^n} Z^2 \prod_{i\le n}\Big(\frac{1}{\sqrt{2\pi}}\exp(-\frac{Y_i^2}{2})dY_i\Big) \\ &= \int_{\mathbb{R}^n}\Bigg(\mathbb{E}_{\theta\sim\pi}\Big[\prod_{i\le n}\exp\Big(Y_i(X\theta)_i - \frac{1}{2}(X\theta)_i^2\Big)\Big]\Bigg)^2 \prod_{i\le n}\Big(\frac{1}{\sqrt{2\pi}}\exp(-\frac{Y_i^2}{2})dY_i\Big) \\ &= \int_{\mathbb{R}^n}\Bigg(\mathbb{E}_{\theta\sim\pi}\Big[\exp(-\frac{1}{2}\|X\theta\|_2^2)\prod_{i\le n}\exp\Big(Y_i(X\theta)_i\Big)\Big]\Bigg)^2 \prod_{i\le n}\Big(\frac{1}{\sqrt{2\pi}}\exp(-\frac{Y_i^2}{2})dY_i\Big) \end{aligned}$$

Recalling (25) we have

$$\begin{aligned} E_0[Z^2] &= \int_{\mathbb{R}^n} \mathbb{E}_{\pi^2}\Big[\exp(-\frac{1}{2}\|X\theta\|_2^2)\exp(-\frac{1}{2}\|X\theta'\|_2^2)\prod_{i\le n}\frac{1}{\sqrt{2\pi}}\exp\Big(Y_i(X(\theta+\theta')_i) - \frac{Y_i^2}{2}\Big)\Big]dy_1\dots dy_n \\ &= \mathbb{E}_{\pi^2}\Big[\exp(-\frac{1}{2}\|X\theta\|_2^2)\exp(-\frac{1}{2}\|X\theta'\|_2^2) \\ &\qquad \prod_{i\le n}\int_{Y_i}\frac{1}{\sqrt{2\pi}}\exp\Big(-\frac{1}{2}(Y_i-(X(\theta+\theta'))_i)^2\Big)dY_i \exp\Big(\frac{1}{2}(X(\theta+\theta'))_i^2\Big)\Big] \\ &= \mathbb{E}_{\pi^2}\Big[\exp\big(\frac{1}{2}\|X(\theta+\theta')\|_2^2 - \frac{1}{2}\|X\theta\|_2^2 - \frac{1}{2}\|X\theta'\|_2^2\big)\Big] \end{aligned}$$

And

$$\begin{aligned} E^X E_0[Z^2] &= E^X\mathbb{E}_{\pi^2}[\exp(\frac{1}{2}\|X(\theta+\theta')\|_2^2 - \frac{1}{2}\|X\theta\|_2^2 - \frac{1}{2}\|X\theta'\|_2^2)] \\ &= \mathbb{E}_{\pi^2}[\exp(\frac{n}{2}(\|\theta+\theta'\|_2^2 - \|\theta\|_2^2 - \|\theta'\|_2^2))\cdot E^X\exp(\frac{1}{2}(Z_1 - Z_2 - Z_3))] \end{aligned} \tag{28}$$

where $Z_1 = \|X(\theta+\theta')\|_2^2 - n\|\theta+\theta'\|_2^2$, $Z_2 = \|X\theta\|_2^2 - n\|\theta\|_2^2$, and $Z_3 = \|X\theta'\|_2^2 - n\|\theta'\|_2^2$.
Using Cauchy–Schwarz Inequality we have:

$$E^X\exp(\frac{1}{2}(Z_1 - Z_2 - Z_3) \le (E^X\exp(Z_1))^{\frac{1}{2}}(E^X\exp(2Z_2))^{\frac{1}{4}}(E^X\exp(2Z_3))^{\frac{1}{4}} \tag{29}$$

Now we evaluate $E^X[\exp(Z_1)]$. Recall that $X_{ij} \overset{\text{i.i.d.}}{\sim} N(0,1)$, so $\|X(\theta+\theta')\|_2^2$ follows a chi-squared distribution:

$$\|X(\theta+\theta')\|_2^2 \sim \|\theta+\theta'\|_2^2 \sum_{i=1}^n g_i^2$$

where $g_i \overset{\text{i.i.d.}}{\sim} N(0,1)$. So

$$Z_1 \sim \|\theta+\theta'\|_2^2 \sum_{i=1}^n (g_i^2 - 1)$$

Then apply Theorem 3.1.9 in [6]: In this case we have $\tau_i = 1$ for any i, so $A = I, \|A\| = 1, \|A\|_{HS}^2 = n$. For $\|\theta+\theta'\|_2^2 < \frac{1}{2}$:

$$E^X \exp(Z_1) \leq \exp\Big(\varphi_{2n,2}(\|\theta+\theta'\|_2^2)\Big) = \exp\Big(\frac{n\|\theta+\theta'\|_2^4}{1-2\|\theta+\theta'\|_2^2}\Big)$$

Very similarly, for $\|\theta\|_2^2, \|\theta'\|_2^2 < \frac{1}{4}$:

$$E^X[\exp(2Z_2)] \leq \exp\Big(\frac{2n\|\theta\|_2^4}{1-4\|\theta\|_2^2}\Big), \quad E^X[\exp(2Z_3)] \leq \exp\Big(\frac{2n\|\theta'\|_2^4}{1-4\|\theta'\|_2^2}\Big)$$

The above shows that for any $\theta \sim \pi$ s.t. $\|\theta\|_2 = o_P(n^{-\frac{1}{4}})$,

$$E^X \exp\Big(\frac{1}{2}(Z_1 - Z_2 - Z_3)\Big) \leq 1 + o(1)$$

And that

$$\begin{aligned}
\mathbb{E}_0[Z^2] = E^X E_0[Z^2] &\leq (1+o(1))\mathbb{E}_{\pi^2}[\exp(\frac{n}{2}(\|\theta+\theta'\|_2^2 - \|\theta\|_2^2 - \|\theta'\|_2^2))] \\
&= (1+o(1))\mathbb{E}_{\pi^2} \exp n \sum_{i=1}^p (\theta \cdot \theta') \qquad (30)
\end{aligned}$$

For $\mathbb{E}_{\pi^2}\big[\exp n \sum_{i=1}^p (\theta\cdot\theta')\big] \leq 1$,

$$\mathbb{E}_0(Z-1)^2 = \mathbb{E}[Z^2] - 1 \leq o(1)$$

And by (26) we have the testing error is no less than 1 for all $\psi$.

□

**Lemma 4.1** (Variance of the $U$-statistic and the linear term)**.** *Under model* (1) *with* $X_{ij} \overset{i.i.d.}{\sim} N(0,1)$*,* $\varepsilon_i \overset{i.i.d.}{\sim} N(0,1)$ *independent of* $X$*, let* $\hat\theta$ *be given as in Theorem 2.1,* $U_n^{-\alpha}(\cdot)$ *and* $L_n^{-\alpha}(\cdot)$ *be defined as in* (20)*. Then for all* $0 < b < \infty, \theta \in \{\theta : \sum_{j=1}^p \theta_j^2 \leq b^2\}$ *and* $\alpha \geq 1/4$ *there exists universal constant C such that*

$$Var_\theta(U_n^{-\alpha}(\theta)) \leq C(1+b^4)\frac{\log p}{n^2}$$

*and*

$$Var_\theta(L_n^{-\alpha}(\hat\theta)) \leq C(1+b^2)\frac{(\log p)^2}{n^2}$$

*with high probability under the law of* $\hat\theta$*.*

*Proof of Lemma 4.1.* Without loss of generality we assume the sample size is $2n$ and the first half of the sample is used to generate the estimation $\hat\theta$. Let $(Y, X)$ denote the second half of the sample, and we take expectations w.r.t. only $(Y, X)$. Noting that $\mathbb{E}_\theta[Y_i X_{ij}] = \theta_j$ and that $\hat\theta$ is independent of $(Y, X)$, for $L_n^{-\alpha}(\cdot)$ we have $\mathbb{E}_\theta[L_n^{-\alpha}(\hat\theta)] = 0$ and $Var_\theta[L_n^{-\alpha}(\hat\theta)] = \mathbb{E}_\theta[L_n^{-\alpha}(\hat\theta)^2]$

We can write

$$(L_n^{-\alpha}(\hat\theta))^2 = \frac{1}{n^2}\sum_{i=1}^n\sum_{j=1}^p\sum_{l=1}^n\sum_{m=1}^p (\theta_j - Y_iX_{ij})(\theta_j - \hat\theta_j)\cdot j^{-2\alpha}(\theta_m - Y_lX_{lm})(\theta_m - \hat\theta_m)\cdot m^{-2\alpha}$$

For $j \neq m$ or $i \neq l$, $(\theta_j - Y_iX_{ij})$ and $(\theta_m - Y_lX_{lm})$ are both zero-mean and mutually independent. So

$$\mathbb{E}_\theta[L_n^{-\alpha}(\hat\theta)^2)] = \frac{1}{n^2}\sum_{j=1}^p (\theta_j - \hat\theta_j)^2 \cdot j^{-4\alpha}\cdot \mathbb{E}_\theta\Big[\sum_{i=1}^n (\theta_j - Y_iX_{ij})^2\Big]$$

In our model $Y_i = \sum_{h=1}^p X_{ih}\theta_h + \varepsilon_i$ and we assume $\|\theta\|_2^2 \coloneqq \sum_{j=1}^p \theta_j^2 \le b^2$ for some known $b$. So for any $1 \le j \le p$:

$$\begin{aligned}
&\mathbb{E}_\theta\Big[\sum_{i=1}^n (\theta_j - Y_i X_{ij})^2\Big] \\
&= \mathbb{E}_\theta\Big[\sum_{i=1}^n \Big\{\theta_j^2 - 2X_{ij}\theta_j\Big(\sum_{h=1}^p X_{ih}\theta_h + \varepsilon_i\Big) + \Big(\sum_{h=1}^p X_{ih}\theta_h + \varepsilon_i\Big)^2\Big\}\Big] \\
&= \mathbb{E}_\theta\Big[\sum_{i=1}^n \Big\{\theta_j^2 - 2X_{ij}\theta_j\varepsilon_i - 2X_{ij}\theta_j\sum_{h=1}^p X_{ih}\theta_h + \varepsilon_i^2 + 2\varepsilon_i\sum_{h=1}^p X_{ih}\theta_h + \Big(\sum_{h=1}^p X_{ih}\theta_h\Big)^2\Big\}\Big] \\
&= \sum_{i=1}^n \Big\{\theta_j^2 - 2\theta_j^2 + 1 + \sum_{h=1}^p \theta_h^2\Big\} \le (1+b^2)n
\end{aligned}$$

By our choice of $\hat{\theta}$, with high probability we have $\|\hat{\theta} - \theta\|_\infty \le \overline{C}\frac{\log p}{n}$ for some constant $\overline{C}$. So for $\alpha \ge 1/4$ there is, with high $\mathbb{P}_\theta$-probability,

$$\begin{aligned}
\mathbb{E}_\theta[L_n^{-\alpha}(\hat{\theta})^2)] &= \frac{1}{n^2}\sum_{j=1}^p (\theta_j - \hat{\theta}_j)^2 \cdot j^{-4\alpha} \cdot \mathbb{E}_\theta\Big[\sum_{i=1}^n (\theta_j - Y_i X_{ij})^2\Big] \\
&\le \frac{1}{n^2}\cdot(1+b^2)n \cdot \sum_{j=1}^p j^{-4\alpha}(\theta_j - \hat{\theta}_j)^2 \\
&\le (1+b^2)C'\frac{(\log p)^2}{n^2}
\end{aligned}$$

where $C'$ takes value of, for example, $2\overline{C}^2$. The second to last line also leads to

$$Var_\theta(L_n^{-\alpha}(\tilde{\theta})) \le \frac{1+b^2}{n}\|\tilde{\theta} - \theta\|^2_{\ell_2^{-\alpha}} \tag{31}$$

for any estimator $\tilde{\theta}$, which might be of separate interest.

For $U_n^{-\alpha}(\theta)$, notice that for any $l \ne m$

$$\mathbb{E}_\theta[(Y_l X_{lj} - \theta_j)(Y_m X_{mj} - \theta_j)] = 0.$$

So $\mathbb{E}_\theta[U_n^{-\alpha}(\theta)] = 0$ and $Var_\theta[U_n^{-\alpha}(\theta)] = \mathbb{E}_\theta[U_n^{-\alpha}(\theta)^2]$.

Denote by $v(l,j) \coloneqq Y_l X_{lj} - \theta_j$,

$$\mathbb{E}_\theta[U_n^{-\alpha}(\theta)^2] = \frac{4}{n^2(n-1)^2}\sum_{j=1}^p\sum_{l=2}^n\sum_{m=1}^{l-1}\sum_{r=1}^p\sum_{s=2}^n\sum_{t=1}^{s-1} j^{-2\alpha}r^{-2\alpha}\mathbb{E}_\theta[v(l,j)v(m,j)v(s,r)v(r,t)]$$

As $Y_l X_{lj} = (\sum_{u=1}^p X_{lu}\theta_u + \varepsilon_l)X_{lj}$, When $j \ne r$ notice that in the sum only when $u = j$ or $u = r$ the coefficients are potentially correlated to coefficients in $v(m,j)$, $v(s,r)$ and $v(r,t)$. So let

$$\bar{v}_{(i,h)}(l,j) \coloneqq (X_{li}\theta_i + X_{lh}\theta_h + \varepsilon_l)X_{lj} - \theta_j$$

and

$$\mathbb{E}_\theta[v(l,j)v(m,j)v(s,r)v(r,t)] = \mathbb{E}_\theta[\bar{v}_{(j,r)}(l,j)\bar{v}_{(j,r)}(m,j)\bar{v}_{(j,r)}(s,r)\bar{v}_{(j,r)}(t,r)]$$

CASE 1: $j \ne r$ and $l, m, s, t$ take four different values.

All $\bar{v}_{(j,r)}(\cdot,\cdot)$ terms inside the expectation are mutually independent and have zero mean. So

$$\mathbb{E}_\theta[\bar{v}_{(j,r)}(l,j)\bar{v}_{(j,r)}(m,j)\bar{v}_{(j,r)}(s,r)\bar{v}_{(j,r)}(t,r)] = 0$$

CASE 2: $j \ne r, l > m = s > t$ or $s > t = l > m$

When $l > m = s > t$, $\bar{v}_{(j,r)}(l,j)$ has zero mean and is independent from all other $\bar{v}_{(j,r)}(\cdot,\cdot)$ terms inside the expectation. So

$$\mathbb{E}_\theta[\bar{v}_{(j,r)}(l,j)\bar{v}_{(j,r)}(m,j)\bar{v}_{(j,r)}(s,r)\bar{v}_{(j,r)}(t,r)] = 0$$

Similarly, when $s > t = l > m$, the fact that $\bar{v}_{(j,r)}(s,r)$ has zero mean and is independent from all other $\bar{v}_{(j,r)}(\cdot,\cdot)$ terms leads to the expectation being zero.

CASE 3: $j \neq r, m = t < s = l$

Notice $\bar{v}_{(j,r)}(m,j)\bar{v}_{(j,r)}(t,r)$ and $\bar{v}_{(j,r)}(l,j)\bar{v}_{(j,r)}(s,r)$ are mutually independent.

Since $s = l$, and $\mathbb{E}_\theta[\varepsilon_l^2 X_{lj} X_{lr}] = 0$, we have

$$\begin{aligned}
&\mathbb{E}_\theta[\bar{v}_{(j,r)}(l,j)\bar{v}_{(j,r)}(s,r)] \\
&= \mathbb{E}_\theta[(X_{lj}^2\theta_j + X_{lr}X_{lj}\theta_r + \varepsilon_l X_{lj} - \theta_j)(X_{lr}^2\theta_r + X_{lr}X_{lj}\theta_j + \varepsilon_l X_{lr} - \theta_r)] \\
&= \mathbb{E}_\theta[X_{lj}^2 X_{lr}^2\theta_j\theta_r + X_{lj}^3 X_{lr}\theta_j^2 - X_{lj}^2\theta_j\theta_r + X_{lr}^3 X_{lj}\theta_r^2 \\
&\qquad + X_{lj}^2 X_{lr}^2\theta_j\theta_r - X_{lr}X_{lj}\theta_r^2 - X_{lr}^2\theta_j\theta_r - X_{lr}X_{lj}\theta_j^2 + \theta_r\theta_j] \\
&= \theta_j\theta_r
\end{aligned}$$

The same goes for $\mathbb{E}_\theta[\bar{v}_{(j,r)}(l,j)\bar{v}_{(j,r)}(s,r)]$ so

$$\mathbb{E}_\theta[\bar{v}_{(j,r)}(l,j)\bar{v}_{(j,r)}(m,j)\bar{v}_{(j,r)}(s,r)\bar{v}_{(j,r)}(t,r)] = \theta_j^2\theta_r^2$$

There are $\frac{n(n-1)}{2}$ combinations of $m = t < s = l$, so the contribution of the corresponding terms sums up to

$$\frac{2}{n(n-1)}\Big(\sum_{r=1}^{p}\sum_{j=1}^{p}\theta_j^2\theta_r^2 \cdot j^{-2\alpha}r^{-2\alpha}\Big) \leq \frac{2b^4}{n(n-1)}$$

When $j = r$, in $v(l,j)$ all but three terms have zero mean and are independent from all other terms in other $v$. Let $\bar{v}_{(i)}(l,j) := X_{li}^2\theta_i + \varepsilon_l X_{lj} - \theta_j$, then

$$\mathbb{E}_\theta[v(l,j)v(m,j)v(s,r)v(r,t)] = \mathbb{E}_\theta[\bar{v}_{(j)}(l,j)\bar{v}_{(j)}(m,j)\bar{v}_{(j)}(s,r)\bar{v}_{(j)}(t,r)]$$

CASE 4: $j = r$ and $l, m, s, t$ take four different values.

All $\bar{v}_{(j)}$ terms inside the expectation are mutually independent and have zero mean. So

$$\mathbb{E}_\theta[\bar{v}_{(j)}(l,j)\bar{v}_{(j)}(m,j)\bar{v}_{(j)}(s,r)\bar{v}_{(j)}(t,r)] = 0$$

CASE 5: $j = r, l > m = s > t$ or $s > t = l > m$.

Same as in CASE 2, evaluating $\bar{v}_{(j)}(l,j)$ when $l > m = s > t$, or $\bar{v}_{(j)}(s,j)$ when $s > t = l > m$ leads to

$$\mathbb{E}_\theta[\bar{v}_{(j)}(l,j)\bar{v}_{(j)}(m,j)\bar{v}_{(j)}(s,r)\bar{v}_{(j)}(t,r)] = 0$$

CASE 6: $j = r, m = t < s = l$

Notice $\bar{v}_{(j)}(m,j)\bar{v}_{(j)}(t,r)$ and $\bar{v}_{(j)}(l,j)\bar{v}_{(j)}(s,r)$ are mutually independent.

$$\begin{aligned}
&\mathbb{E}_\theta[\bar{v}_{(j)}(m,j)\bar{v}_{(j)}(t,r)] \\
&= \mathbb{E}_\theta[X_{mj}^4\theta_j^2 - X_{mj}^2\theta_j^2 + \varepsilon_m^2 X_{mj}^2 - X_{mj}^2\theta_j^2 + \theta_j^2] \\
&= 1
\end{aligned}$$

We have $\mathbb{E}_\theta[\bar{v}_{(j)}(l,j)\bar{v}_{(j)}(s,r)] = 1$ by the same procedure.

Again for each $j = r$ there are $\frac{n(n-1)}{2}$ combinations of $m = t < s = l$. So these terms contribute

$$\frac{4}{n^2(n-1)^2}\frac{n(n-1)}{2}\sum_{j=1}^{p} j^{-4\alpha} \leq 2\frac{1+\log p}{n(n-1)}$$

Combining the cases and constants gives $\mathbb{E}[U_n^{-\alpha}(\theta)^2)] \leq 3(1+b^4)\frac{\log p}{n^2}$. For $\alpha > \frac{1}{4}$, the $\log p$ can be improved to a constant depending only on $\alpha$.

□

# A Concentration inequalities related to Gaussian polynomials

Let $g_i, g_i', i \in \mathbb{N}$ be independent $N(0,1)$ random variables. In our context we often need to bound the tail of the Gaussian polynomials $Q(n) := \sum_{i=1}^{n} g_i^2$ and $Z(n) := \sum_{i=1}^{n} g_i \cdot g_i'$. These random variables follow sub-exponential distributions for which concentration inequalities can be derived from Bernstein's inequality (see for example Theorem 3.1.8 in [6]). The two lemmas below give the tail bounds we need for this paper.

**Lemma A.1.** *Let $g_i$, $i \in \mathbb{N}$ be independent $N(0,1)$ random variables and $Q(n) := \sum_{i=1}^{n} g_i^2$. For any $t > 0$*

$$P(|Q(n) - n| \geq \sqrt{nt} + t) \leq C \exp(-C't)$$

*for some $C, C' > 0$.*

The above lemma can be shown by using Theorem 3.1.9 in [6] and taking the matrix $A$ equal to the Identity. The next lemma can be shown either through Bernstein's inequality, or by seeing that $Z(n) = \frac{1}{4}\sum_{i=1}^{n}((g_i + g_i')^2 - (g_i - g_i')^2)$ and using Lemma A.1.

**Lemma A.2.** *Let $g_i, g_i' \overset{i.i.d.}{\sim} N(0,1)$, $i = 1, 2, \ldots, n$ and $Z(n) := \sum_{i=1}^{n} g_i \cdot g_i'$. For any $t > 0$*

$$P\Big(|Z(n)| \geq \sqrt{nt} + t\Big) \leq C \exp(-C't)$$

*for some $C, C' > 0$.*

As squared-Gaussians are sub-exponential, we have a bound for the maximum of $p \simeq n$ independent $Z(n)$'s:

**Lemma A.3.** *Let $Z(n) := \sum_{i=1}^{n} g_i \cdot g_i'$ and $Z_j(n)$, $j \in \mathbb{N}$ be independent draws of $Z(n)$, then assuming that $p \leq C_4 n$ for some constant $C_4$, there exists constant $C_5$ depending on $C_4$ such that*

$$\mathbb{E}[\max_{j \leq p} Z_j(n)] \leq C_5 \sqrt{n \log p}$$

*and for $n \to \infty$*

$$\max_{j \leq p}\{|Z_j(n)|\} = O_P(\sqrt{n \log p})$$

*Proof.* By definition (3.21) and part (b) of Theorem 3.1.10 in [6], and realizing that $Z(n)$ is sub-exponential with $\nu = 2n$ and $c = 2$, we have

$$\begin{aligned}
\mathbb{E}[\max_{j \leq p} Z_j(n)] &\leq \sqrt{4n \log p} + 2 \log p \\
&\leq 2\left(1 + \sqrt{\frac{\log C_4 n}{n}}\right)\sqrt{n \log p} \\
&\leq 2\left(2 + \sqrt{\log C_4}\right)\sqrt{n \log p}
\end{aligned}$$

which gives the first part of the above lemma.

The second half follows by Markov's inequality. □

We have a similar result for concentration of Gram Matrix $X^T X/n$ around its mean:

**Lemma A.4.** *Using same assumptions 2.1 on $X$ and $\theta$ as rest of this paper, for any $1 \leq j \leq p$ and $C < 1$*

$$\mathbb{P}_\theta\left(\Big(\frac{X^T X \theta}{n}\Big)_j - E_\theta\Big[\Big(\frac{X^T X \theta}{n}\Big)_j\Big] > \frac{C-1}{2}\sqrt{\frac{\log p}{n}}\right) \geq 1 - C' \exp(-\log p)$$

*for some $C' > 0$ depending only on $C, b$.*

*Proof.* We have

$$(X^T X\theta)_j = \sum_{l=1}^{p} \Big( \sum_{i=1}^{n} X_{il} X_{ij} \Big) \theta_l$$

and

$$\mathbb{E}[(X^T X\theta)_j] = \theta_j$$

So

$$\Big(\frac{X^T X\theta}{n}\Big)_j - E_\theta\Big[\Big(\frac{X^T X\theta}{n}\Big)_j\Big] = \Big(\frac{1}{n}\sum_{i=1}^{n} X_{ij}^2 \theta_j - \theta_j\Big) + \Big(\frac{1}{n}\sum_{l\neq j}^{p}(\sum_{i=1}^{n} X_{il} X_{ij})\theta_l\Big)$$

Concentration inequality for first part can be established by Lemma A.1. For the second half we have:

$$\frac{1}{n}\sum_{l\neq j}^{p}(\sum_{i=1}^{n} X_{il} X_{ij})\theta_l \leq \frac{1}{n}\|\theta\|_1 \cdot \max_l(\sum_{i=1}^{n} X_{il} X_{ij}),$$

tail probability of which can be bounded by Lemma A.3 since $\|\theta\|_1$ is bounded by assumption B in 2.1. □